\documentclass[12pt]{article}
\usepackage{amsfonts,graphicx,color,amsmath,amssymb}

\newcommand{\partm}{\mathop{\mathrm{pm}}}
\newcommand{\parti}{\mathop{\mathrm{pi}}}
\newcommand{\Aut}{\mathop{\mathrm{Aut}}}
\newcommand{\lcm}{\mathop{\mathrm{lcm}}}

\newtheorem{theorem}{Theorem}[section]
\newtheorem{prop}[theorem]{Proposition}
\newtheorem{cor}[theorem]{Corollary}
\newtheorem{lemma}[theorem]{Lemma}
\newtheorem{unnumber}{}

\newtheorem{prepf}[unnumber]{Proof}
\newenvironment{pf}{\prepf\rm}{\endprepf}
\newcommand{\qed}{\hfill$\Box$}

\begin{document}

\title{A graph partition problem}
\author{Peter J. Cameron\footnote{Mathematical Institute, North Haugh, St. Andrews KY16 9SS, UK; {\tt pjc@mcs.st-andrews.ac.uk}} \, and Sebastian M. Cioab\u{a}\footnote{Department of Mathematical Sciences, University of Delaware, Newark, DE 19716-2553, USA; {\tt cioaba@math.udel.edu}}}
\date{July 30, 2014}
\maketitle

\begin{abstract}
Given a graph $G$ on $n$ vertices, for which $m$ is it possible to partition
the edge set of the $m$-fold complete graph $mK_n$ into copies of $G$?
We show that there is an integer $m_0$, which we call the \emph{partition
modulus of $G$}, such that the set $M(G)$ of values of $m$ for which
such a partition exists consists of all but finitely many multiples of $m_0$.
Trivial divisibility conditions derived from $G$ give an integer $m_1$ which
divides $m_0$; we call the quotient $m_0/m_1$ the \emph{partition index of $G$}. It seems that most graphs $G$ have partition index equal to $1$, but we
give two infinite families of graphs for which this is not true. We also
compute $M(G)$ for various graphs, and outline some connections between
our problem and the existence of designs of various types.
\end{abstract}

\section{Introduction}

The problem of interest in this paper is the following:

Given a graph $G$ on $n$ vertices, is it possible to partition the edge set of the complete graph $K_n$ into isomorphic copies of $G$ ? If this is not possible, then we are interested in determining the set of integers $m$ such that the edge set of the $m$-fold complete graph $mK_n$ can be partitioned into copies of $G$.

We will see that this seemingly simple problem has connections to algebra, combinatorics and geometry among others. Important open problems such as the
existence of finite projective planes are equivalent to edge partition
problems into certain specified graphs.

\begin{center}
\begin{figure}[htbp]
\includegraphics[scale=0.5]{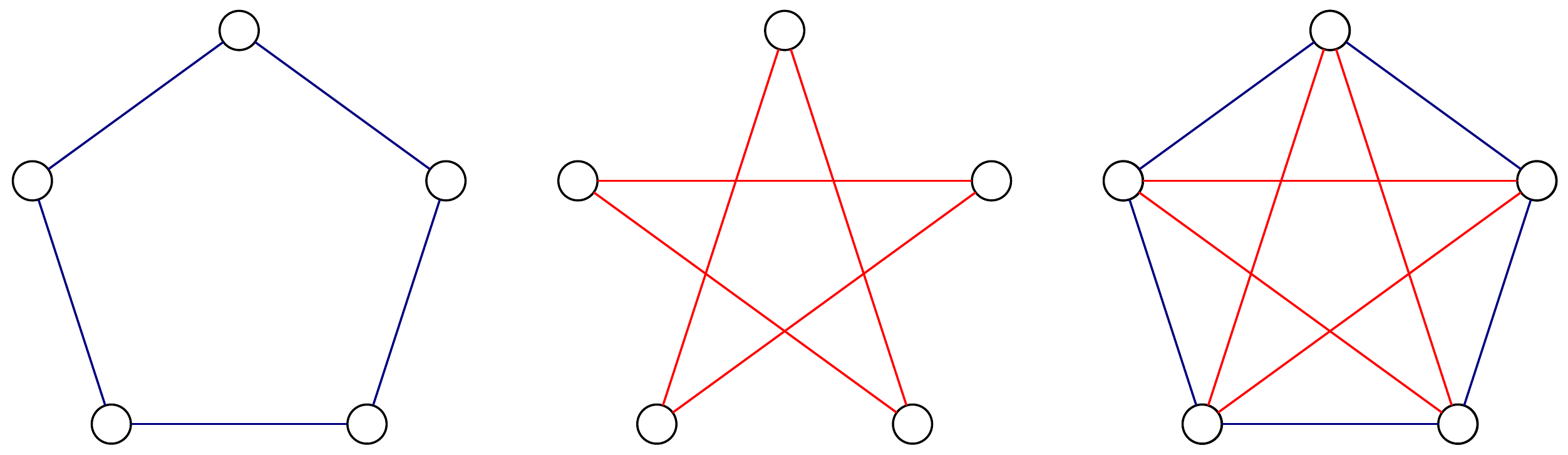}
\caption{\label{f:twoC5}Two edge disjoint cycles on $5$ vertices on the same vertex set.}
\end{figure}
\end{center}

Historically, perhaps the first instance of this type problem goes back to Walecki \cite{lucas} who showed that the edge-set of the complete graph $K_n$ can be partitioned into copies of the cycle $C_n$ when $n\geq 3$ is odd and the edge-set of the complete graph $K_n$ minus a perfect matching can be partitioned into copies of $C_n$ when $n\geq 4$ is even. See Figure~\ref{f:twoC5} for a partition of $K_5$ into edge disjoint copies of $C_5$. This result shows the impossibility of decomposing the edge of $K_n$ into copies of $C_n$, but also that the edge set of $2K_n$ can be partitioned into copies of $C_n$, when $n\geq 4$ is even. 

\begin{center}
\begin{figure}[htbp]
\includegraphics[scale=0.6]{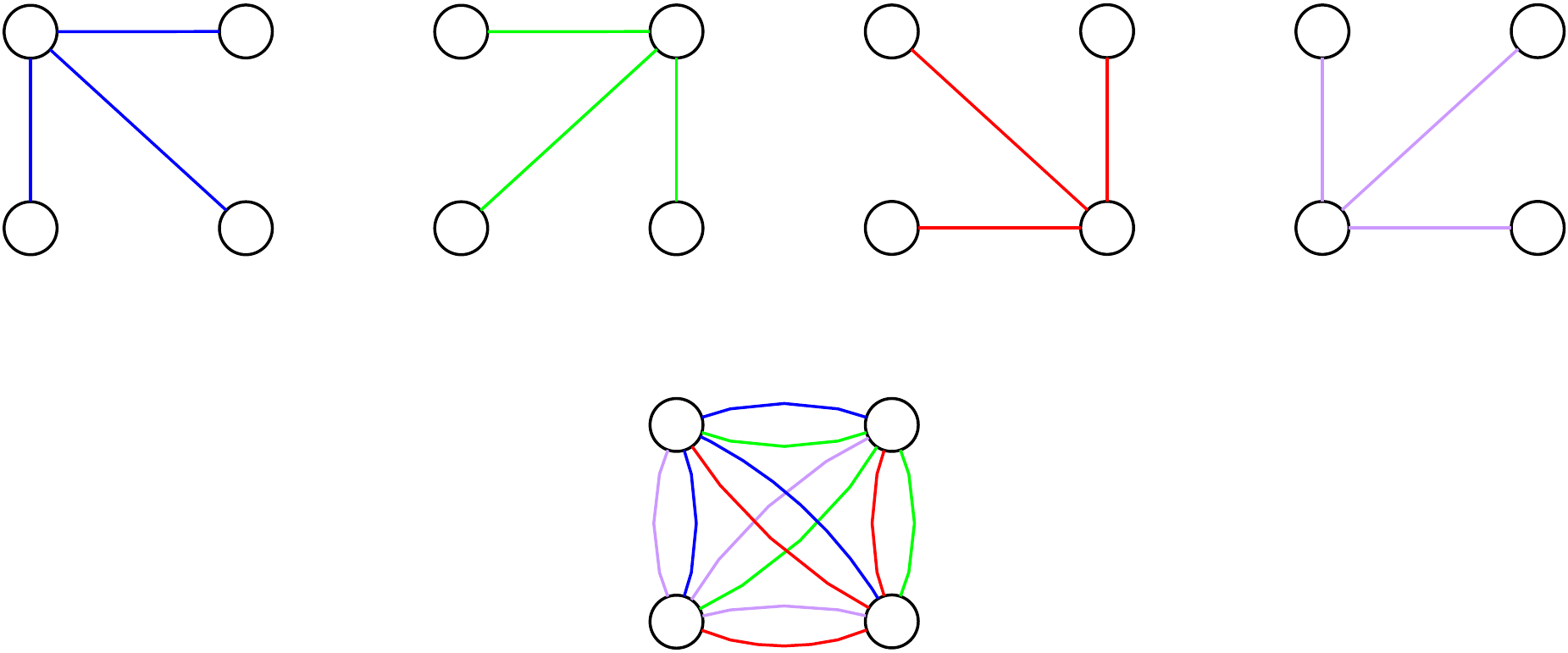}
\caption{\label{f:twoK4}Four stars $K_{1,3}$ decomposing $2K_4$.}
\end{figure}
\end{center}

Figure~\ref{f:twoK4} shows another instance of this problem where $K_4$ cannot be decomposed into edge disjoint copies of $K_{1,3}$, but $2K_4$ can be partitioned into $4$ $K_{1,3}$s.

\begin{center}
\begin{figure}[htbp]
\includegraphics[scale=0.7]{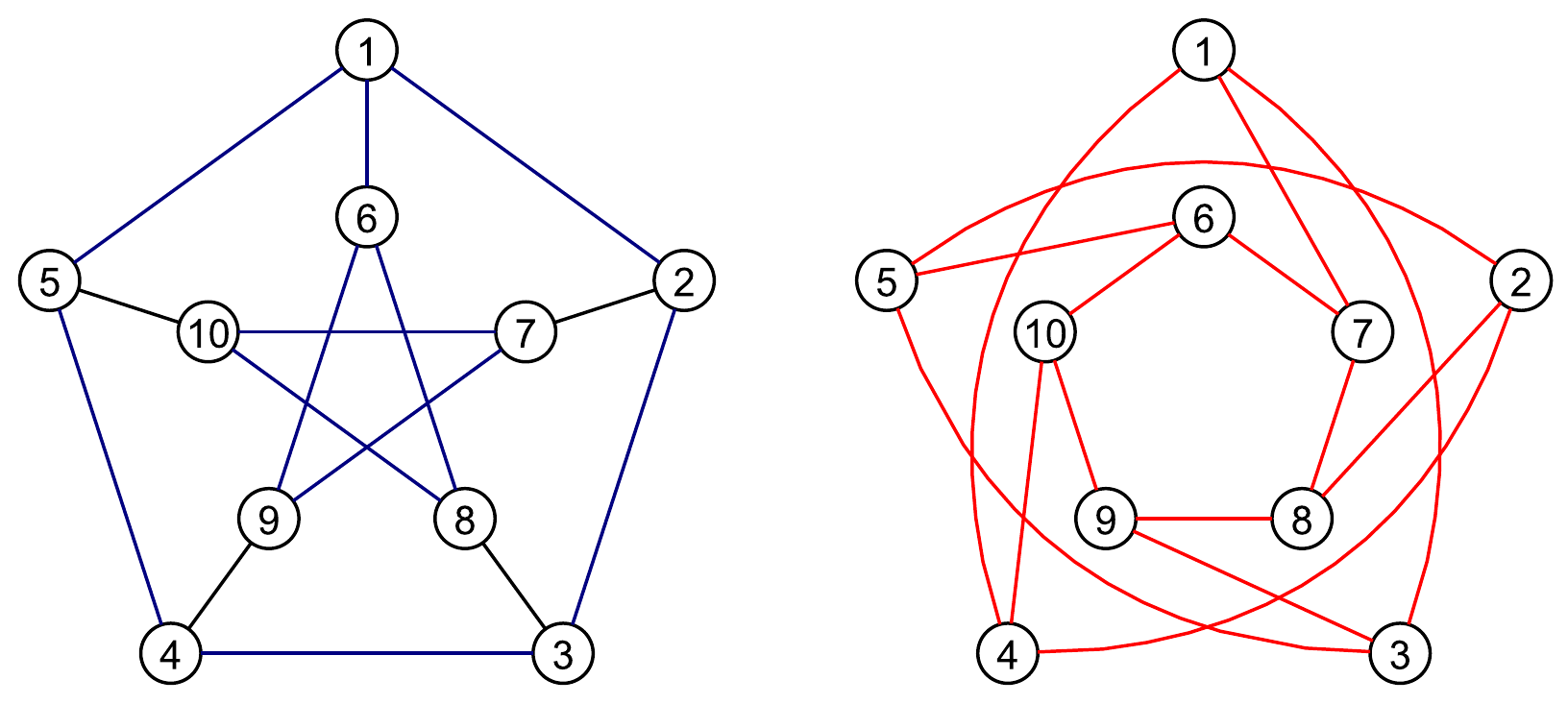}
\caption{\label{f:twopeters}Two edge disjoint Petersen graphs on the same vertex set.}
\end{figure}
\end{center}

Figure~\ref{f:twopeters} shows that it is possible to find two edge-disjoint
copies of the Petersen graph on $10$ vertices. But can we find three edge-disjoint copies, that is, a partition of the edges of the complete graph on $10$ vertices into three Petersen graphs?

This problem was proposed by Allen Schwenk in the \textit{American Math Monthly} in 1983. Elegant negative solutions by Schwenk and O.~P.~Lossers appeared in the same journal in 1987~\cite{schwenk}.
This result is described in many books on algebraic graph theory including
Godsil and Royle \cite[Section 9.2]{GR}, and Brouwer and Haemers
\cite[Section 1.5.1]{BH}. Schwenk's argument plus some simple combinatorial ideas can be used to show that whenever we can arrange two edge disjoint Petersen graphs on the same vertex set, then the
complement of their union must be the bipartite cubic graph on $10$ vertices that is the bipartite complement of $C_{10}$. In Figure \ref{f:twopeters}, every missing edge goes from the set of $5$ outer vertices to the set of $5$ inner vertices, so the complement of the union of the two Petersen graphs is visibly bipartite (see \cite{2Pet} for a proof).

A generalization was posed by 
Rowlinson~\cite{rowlinson}, and further variants have also been studied.
\v{S}iagiov\'{a} and Meszka \cite{sm} obtained a packing of five 
Hoffman--Singleton graphs in the complete graph $K_{50}$. At present time,
it is not known if it is possible to decompose $K_{50}$ into seven 
Hoffman--Singleton graphs. Van~Dam \cite{vandam1} showed that if the
edge-set of the complete graph of order $n$ can be partitioned into three
(not necessarily isomorphic) strongly regular graphs of order $n$, then this
decomposition forms an amorphic association scheme (see also van Dam and
Muzychuk \cite{vandammuzychuk}). 

At the Durham Symposium on Graph Theory and Interactions in 2013, the
authors amused themselves by showing that, for every $m>1$, the $m$-fold
complete graph $mK_{10}$
(with $m$ edges between each pair of vertices) can be partitioned into
$3m$ copies of the Petersen graph. Our purpose in this paper is to extend
this investigation by replacing the Petersen graph by an arbitrary graph.

In fact, a stronger result for the Petersen graph was found by
Adams, Bryant and Khodkar \cite{abk}; these authors allow $n$ to
be arbitrary, in other words, they allow adding arbitrary many isolated vertices
to the Petersen graph. We do not consider this more general problem.

\section{Partition Modulus}

\paragraph{Definition}
For a graph $G$ on $n$ vertices, we let
\[M(G)=\{m:mK_n\hbox{ can be partitioned into copies of }G\}.\]

\paragraph{Example}
As mentioned in the Introduction, one of the earliest results on this concept is that of Walecki
\cite{lucas}, according to which the complete graph $K_n$ can be partitioned
into $(n-1)/2$ Hamiltonian cycles if (and only if) $n$ is odd; if $n$ is
even, then $2K_n$ can be partitioned into Hamiltonian cycles. If $n$ is even, $2K_n$ can be partitioned into Hamiltonian cycles, but $mK_n$ cannot if $m$ is odd, since $n$ does not divide $mn(n-1)/2$ for
$n$ even and $m$ odd. So
\[
M(C_n)=\begin{cases} \mathbb{N} & \text{if } n \text{ is odd},\\
                    2\mathbb{N} & \text{if } n \text{ is even}.\\\end{cases}
\]

\begin{prop}
For any graph $G$, the set $M(G)$ is non-empty. In fact, if $G$ has $e$
edges and $\Aut(G)$ is its automorphism group, then
\[2(n-2)!e/|\Aut(G)|\in M(G).\]
\label{p:exist}
\end{prop}

\begin{pf}
The graph $G$ has $n!/|\Aut(G)|$ images under the symmetric group $S_n$, since
$S_n$ acts transitively on the set of graphs on vertex set $\{1,\ldots,n\}$
which are isomorphic to $G$, and $\Aut(G)$ is the stabiliser of one of these
graphs. Each of the $n(n-1)/2$ pairs of points is covered equally often by
an edge in one of these images, since $S_n$ is doubly transitive; double
counting gives this number to be
$\Big(n!/|\Aut(G)|\Big)e\Big/\Big(n(n-1)/2\Big)$,
as required.\qed
\end{pf}

Our main result is a description of the set $M(G)$.

\begin{theorem}
For any graph $G$, there is a positive integer $m_0$ and a finite set $F$ of
multiples of $m_0$ such that $M(G)=m_0\mathbb{N}\setminus F$.
\end{theorem}

We call the number $m_0$ the \emph{partition modulus} of $G$, and denote it
by $\partm(G)$.

The theorem follows immediately from a couple of simple lemmas.

\begin{lemma}
The set $M(G)$ is additively closed.
\end{lemma}

\begin{pf}Superimposing partitions of the edges of $aK_n$ and $bK_n$
gives a partition of $(a+b)K_n$.
\qed\end{pf}

\begin{lemma}
An additively closed subset $M$ of $\mathbb{N}$ has the form
$m\mathbb{N}\setminus F$, where $F$ is a finite set of multiples of $m$.
\end{lemma}

\begin{pf}
We have no convenient reference (though \cite{nw} is related), so we sketch
the proof. We let $m=\gcd(M)$. By dividing through by $m$, we obtain a set
with gcd equal to $1$, so it suffices to prove the result in this case.

First we observe that $M$ is finitely generated, that is, there is a finite
subset $K$ such that any element $M$ is a linear combination, with non-negative
integer coefficients, of elements of $K$. Then we proceed by induction on
$|K|$. It is well known that, if $\gcd(a,b)=1$, then all but finitely many
positive integers have the form $xa+yb$ for some $x,y\ge0$. Assume that the
result holds for generating sets smaller than $K$. Take $a\in K$, and let
$b=\gcd(K\setminus\{a\})$. By induction, $K\setminus\{a\}$ generates all but
finitely many multiples of $b$. Also, $\gcd(a,b)=1$, so that the result for
sets of size $2$ finishes the argument.\qed
\end{pf}

We have not tried to get an explicit bound here, since for most graphs the
excluded set $F$ seems to be much smaller than our general argument suggests.

\paragraph{Example} As mentioned in Section \ref{intro}, if $P$ is the Petersen graph, then
$M(P)=\mathbb{N}\setminus\{1\}$, so that the partition modulus of the
Petersen graph is $1$. It suffices to show that $2,3\in M(P)$.

That $2\in M(P)$ follows from a generalisation of Proposition~\ref{p:exist}:

\begin{prop}
Suppose that $G$ has $n$ vertices and $e$ edges, and that there is a
doubly transitive group $H$ of degree $n$ for which $|H:H\cap\Aut(G)|=r$.
Then $2re/n(n-1)\in M(G)$.
\label{p:2tr}
\end{prop}

\begin{pf}
The graph $G$ has $r$ images under $H$, whose $re$ edges cover all pairs
$2re/n(n-1)$ times.\qed
\end{pf}

Now $\Aut(P)\cong S_5$, a subgroup of index $6$ in $S_6$ (which acts as a
$2$-transitive group on the vertex set of $P$). So $6\cdot15/45=2\in M(P)$.

A direct construction shows that $3\in M(P)$. We do this by means of a
$9$-cycle on the
vertex set of $P$, fixing a point $\infty$ and permuting the remaining points
as $(0,1,2,\ldots,8)$. It is clear that the images of the three edges of $P$
containing $\infty$ cover all pairs of the form $\{\infty,x\}$ three times.
For the remaining pairs, we need to choose a drawing of $P-\infty$ on the
vertices $\{0,\ldots,8\}$ in such a way that each of the distances $1,2,3,4$
in the $9$-cycle is represented three times by an edge, since these distances
index the orbits of the cycle on $2$-sets. It takes just a moment by computer
to find dozens of solutions. For example, the edges
\begin{eqnarray*}
&& \{\infty,1\},\{\infty,4\},\{\infty,8\}, \{0,2\},\{0,3\},\{0,4\},\{1,3\},
\{1,5\},\\
&& \{2,5\},\{2,8\},\{3,7\},\{4,6\},\{5,6\},\{6,7\},\{7,8\}
\end{eqnarray*}
have the required properties.

\section{Partition index}

As is common in problems of this kind, there are some divisibility conditions
which are necessary for a partition to exist:

\begin{prop}
Let $G$ have $n$ vertices and $e$ edges. Then every element $m$ of $M(G)$
has the property that $e$ divides $mn(n-1)/2$, and the greatest common divisor
of the vertex degrees of $G$ divides $m(n-1)$.
\label{p:div}
\end{prop}

\begin{pf}
If $l$ copies of $G$ partition $mK_n$, then $mn(n-1)/2=le$, proving the
first assertion; and the $m(n-1)$ edges through a vertex in $mK_n$ are
partitioned by the vertex stars in copies of $G$, from which the second
assertion follows.\qed
\end{pf}

Let $m_1$ be the number for which these divisibility conditions are
equivalent to the assertion that $m_1\mid m$ for all $m\in M(G)$.
Thus,
\[m_1=\lcm\left(\frac{e}{\gcd(e,n(n-1)/2)},\frac{d}{\gcd(d,n-1)}\right),\]
where $d$ is the greatest common divisor of the vertex degrees.

We have that $m_1\mid m_0=\partm(G)$. We define the \emph{partition
index} $\parti(G)$ to be the quotient $m_0/m_1$.

\begin{prop}
Let $G$ have $n$ vertices and $e$ edges, and let $\overline{G}$ denote its
complement. Then $m\in M(G)$ if and only if
$m\Big(n(n-1)/2-e\Big)\Big/e\in M(\overline{G})$.
\end{prop}

\begin{pf}
If $l$ copies of $G$ cover $kK_n$, then $l$ copies of $\overline{G}$ cover
$(l-k)K_n$. Using $l=kn(n-1)/(2e)$ from the preceding Proposition gives
the result.\qed
\end{pf}

\begin{cor}
If $G$ has $n$ vertices and $e$ edges, then
\[\partm(\overline{G})=\left(\frac{n(n-1)/2-e}{e}\right)\partm(G).\]
\label{p:comp}
\end{cor}

The \emph{triangular graph} $T(l)$ is the line graph of $K_l$, that is,
its vertices are the $2$-element subsets of an $l$-set, two vertices
adjacent if they have non-empty intersection.

\paragraph{Example}
Since $T(5)$ (the line graph of $K_5$) is the complement of the Petersen
graph, we have $M(T(5))=2\mathbb{N}\setminus\{2\}$.

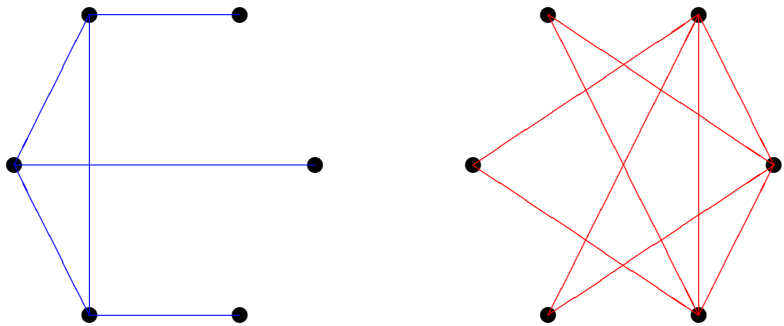
\begin{figure}[htbp]
\begin{center}
\setlength{\unitlength}{2mm}
\begin{picture}(25,20)
\multiput(5,0)(10,0){2}{\circle*{1}}
\multiput(0,10)(20,0){2}{\circle*{1}}
\multiput(5,20)(10,0){2}{\circle*{1}}
\color{blue}
\put(5,0){\line(0,1){20}}
\put(5,0){\line(-1,2){5}}
\put(5,20){\line(-1,-2){5}}
\put(5,0){\line(1,0){10}}
\put(0,10){\line(1,0){20}}
\put(5,20){\line(1,0){10}}
\end{picture}
\qquad
\setlength{\unitlength}{2mm}
\begin{picture}(25,20)
\multiput(5,0)(10,0){2}{\circle*{1}}
\multiput(0,10)(20,0){2}{\circle*{1}}
\multiput(5,20)(10,0){2}{\circle*{1}}
\color{red}
\put(15,0){\line(0,1){20}}
\put(15,0){\line(1,2){5}}
\put(15,20){\line(1,-2){5}}
\put(5,0){\line(3,2){15}}
\put(5,0){\line(1,2){10}}
\put(0,10){\line(3,-2){15}}
\put(0,10){\line(3,2){15}}
\put(5,20){\line(1,-2){10}}
\put(5,20){\line(3,-2){15}}
\end{picture}
\end{center}
\caption{\label{f:twographs}An example}
\end{figure}

\paragraph{Remark}
If the same relation held between the numbers $m_1$ for $G$ and $\overline{G}$
defined earlier as for the partition moduli in Corollary~\ref{p:comp}, then
we would have $\parti(G)=\parti(\overline{G})$. But this
is not true, as we will show using the graphs in Figure~\ref{f:twographs}.

Let $G$ be the graph with $6$ vertices and $6$ edges consisting of a triangle
with a pendant edge at each vertex. Then $G$ has $6$ edges, and the gcd of the
vertex degrees is $1$, so $m_1(G)$ is the least common multiple of
$6/\gcd(6,15)$ and $1/\gcd(1,5)$, that is, $m_1(G)=2$.  However, 
$\overline{G}$ has $9$ edges and all vertex degrees even; so
$m_1(\overline{G})$ is the least common multiple of $9/\gcd(9,15)$ and
$2/\gcd(2,5)$, that is, $m_1(\overline{G})=6$.

We have $4\in M(G)$. This follows from
Proposition~\ref{p:2tr}, since $\Aut(G)$, which is dihedral of order $6$, is a
subgroup of index $10$ in the $2$-transitive group $\mathrm{PSL}(2,5)$. So
$6\in M(\overline{G})$. It follows that $M(\overline{G})=6\mathbb{N}$, and
$\partm(\overline{G})=6$ and $\parti(\overline{G})=1$. But, by
Corollary~\ref{p:comp}, we have $M(G)=4\mathbb{N}$, so that $\partm(G)=4$
and $\parti(G)=2$.

\section{Examples}

In this section, we construct two families of examples of graphs which have
partition index greater than $1$. We are grateful to Mark Walters for the
first of these.

\begin{prop}
Let $G$ be a star $K_{1,n-1}$, with $n>2$. Then $\partm(G)=2$ (and indeed
$M(G)=2\mathbb{N}$); so
$$
\parti(G)=\begin{cases}1 & \text{ if } n \text{ is odd },\cr 2 & \text{ if } n \text{ is even.}\cr\end{cases}
$$
\end{prop}

\begin{pf}
$m_1(G)$ is the least common multiple of $(n-1)/\gcd(n-1,n(n-1)/2)$ (which
is equal to $1$ if $n$ is even, and $2$ if $n$ is odd) and $1/\gcd(1,n-1)=1$.

Suppose that $mK_n$ is covered with copies of the star; let the vertex set
be $\{1,2,\ldots,n\}$, and let $x_i$ be the number of stars with centre at
the vertex $i$. Then the edge $\{i,j\}$ is covered $x_i+x_j$ times; so
$x_i+x_j=m$ for all $i\ne j$. This forces $x_i$ to have a constant value $x$,
and $m=2x$. But we can achieve $m=2$ by taking one star with each possible
centre ($x_i=1$ for all $i$).

So the partition modulus of the star is $2$, and the partition index is
as claimed.\qed
\end{pf}

For the second construction, we observe that, if $n$ is a multiple of $4$,
then the number $n(n-1)/2$ of edges of the complete graph is even, so there
are graphs $G$ and $\overline{G}$ each having $n(n-1)/4$ edges. So the first
term in the lcm for both $m_1(G)$ and $m_1(\overline{G})$ is $1$. Suppose we
arrange that $G$ has all degrees odd; then $\overline{G}$ will have all
degrees even. If $d$ and $\overline{d}$ are the least common multiples of
these degrees, then $d/\gcd(d,n-1)$ is odd, but
$\overline{d}/\gcd(\overline{d},n-1)$ is even (since $n-1$ is odd).
Thus, $m_1(G)$ is odd but $m_1(\overline{G})$ is even.

On the other hand, Proposition~\ref{p:comp} shows that
the partition moduli of the two graphs are equal, and so necessarily even.
Thus, certainly, $\parti(G)>1$.

The smallest example of this construction, for $n=4$, has for $G$ the star
$K_{1,3}$ and for $\overline{G}$ the graph $K_3\cup K_1$. As we saw, the
partition indices of these graphs are $2$ and $1$ respectively. For larger
$n$, the construction gives many examples.

\paragraph{Problem} Is there a simple method of calculating the partition
modulus (and hence the partition index) of a graph $G$? Is it true that
almost all graphs have partition index $1$?

\section{Connection with design theory}

For some special graphs, our partition problem is equivalent to the existence
of certain $2$-designs. A $2$-$(n,k,\lambda)$ design consists of a set
of $n$ points and a collection of $k$-element subsets called blocks, such that
any two points lie in exactly $\lambda$ blocks. The design is resolvable if the
blocks can be partitioned into classes of size $n/k$, each class forming a
partition of the point set.

The existence of $2$-designs has received an enormous amount of study; we
refer to~\cite{bjl} for some results.

Now the following result is clear.

\begin{theorem}
\begin{enumerate}
\item Let $G_1$ be the graph consisting of a $k$-clique and $n-k$ isolated
vertices. Then $m\in M(G_1)$ if and only if there exists a
$2$-$(n,k,m)$ design.
\item Let $k\mid n$, and let $G_2$ be the graph consisting of $n/k$ disjoint
$k$-cliques. Then $m\in M(G_1)$ if and only if there exists a
resolvable $2$-$(n,k,m)$ design.
\end{enumerate}
\end{theorem}

\begin{figure}[htbp]
\begin{center}
\setlength{\unitlength}{2cm}
\begin{picture}(2,2)
\thicklines
\color{black}
\multiput(0,0)(1,0){3}{\line(0,1){2}}
\color{red}
\multiput(0,0)(0,1){3}{\line(1,0){2}}
\color{blue}
\put(2,0){\line(-1,1){2}}
\multiput(0,1)(1,1){2}{\line(1,-1){1}}
\multiput(0,0)(0,1){2}{\line(2,1){2}}
\multiput(0,0)(1,0){2}{\line(1,2){1}}
\color{green}
\put(0,0){\line(1,1){2}}
\multiput(0,1)(1,-1){2}{\line(1,1){1}}
\multiput(0,2)(0,-1){2}{\line(2,-1){2}}
\multiput(0,2)(1,0){2}{\line(1,-2){1}}
\color{black}
\multiput(0,0)(1,0){3}{\circle*{0.1}}
\multiput(0,1)(1,0){3}{\circle*{0.1}}
\multiput(0,2)(1,0){3}{\circle*{0.1}}
\end{picture}
\end{center}
\caption{\label{f:ag23}The affine plane $\mathrm{AG}(2,3)$}
\end{figure}
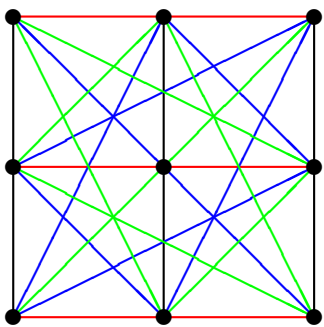

Figure~\ref{f:ag23} shows $K_9$ decomposed into four graphs, each the union
of three disjoint triangles, otherwise known as the affine plane
$\mathrm{AG}(2,3)$. To reduce clutter, we adopt the convention that
a line through three collinear points represents a triangle.

In fact, a more general version of this theorem is true.

\begin{prop}
Let $G$ be a graph on $n$ vertices whose edge-set can be partitioned into
$s$ complete graphs on $k$ vertices. Then a necessary condition for
$\lambda\in M(G)$ is that there exists a $2$-$(n,k,\lambda)$ design.
\label{p:clique}
\end{prop}

This class of graphs includes, for example, the point graphs of partial
geometries and generalized polygons.

This result gives a proof that $1\notin M(T(l))$ (see Theorem \ref{thm:triangular} for another proof): for $T(l)$ is
the edge-disjoint union of $l$ cliques of size $l-1$, and the resulting
$2$-$(l(l-1)/2,l-1,1)$ design would have
\[l(l-1)/2\cdot(l+1)(l-2)/2\Big/(l-1)(l-2)=l(l+1)/4\]
blocks, and so would violate Fisher's inequality (asserting that a $2$-design
has at least as many blocks as points).

In some cases, the converse is true.

\begin{prop}
Let $L_2(q)$ denote the line graph of $K_{q,q}$, the $q\times q$ square lattice
graph. Then $1\in M(L_2(q))$ if and only if $q$ is odd and there exists a
projective plane of order $q$.
\end{prop}

\begin{pf}
We begin by noting that the existence of a projective plane of order $q$ is
equivalent to that of an affine plane of order $q$ (a $2$-$(q^2,q,1)$ design);
such a design is necessarily resolvable, with $q+1$ parallel classes.

Now the necessity of the condition follows from our general results. If the
affine plane exists, then partition the parallel classes into $(q+1)/2$ sets
of size $2$; each set, regarded as a set of $2q$ complete graphs of size $q$,
gives a copy of $L_2(q)$.\qed
\end{pf}

In Figure~\ref{f:ag23}, if we identify red and black, and also blue and green,
we obtain a decomposition of $K_9$ into two copies of $L_2(3)$.

The problem of determining $M(L_2(q))$ in cases not covered by this result,
especially those when the required plane does not exist, is open. For
example, what is $M(L_2(6))$?

It is worth mentioning that the problem is solved
for the unique strongly regular graph with the same parameters as
$L_2(q)$ but not isomorphic to it. This is the $16$-vertex \emph{Shrikhande
graph}. Darryn Bryant \cite{bryant} found five copies of
the Shrikhande graph $S$ that cover the edges of $K_{16}$ twice; so
$M(S)=2\mathbb{N}$, and $\parti(S)=1$.

\section{Triangular graphs}

Recall that the \emph{triangular graph} $T(l)$ is the line graph of $K_l$, that is,
its vertices are the $2$-element subsets of an $l$-set, two vertices
adjacent if they have non-empty intersection.

Now $T(l)$ has $l(l-1)/2$ vertices and valency $2(l-2)$, and has 
$l(l-1)(l-2)/2$ edges. Thus $e/\gcd(e,n(n-1)/2) = 4$, $2$ or $1$ according
as the power of $2$ dividing $l+1$ is $1$, $2$ or at least $4$. Also,
$d=l-2$, so $d/\gcd(d,n-1)$ is $1$ or $2$ according as $l+1$ is even or odd.
So
$$
m_1(T(l))=\begin{cases}1 & \text{ if } l\equiv3\pmod{4},\cr
              2 & \text{ if } l\equiv1\pmod{4},\cr
              4 & \text{ if } l \text{ is even.}\cr\end{cases}$$
We conjecture that these are also the partition moduli of the triangular
graphs, so that the partition indices are all $1$. We saw this already for
$T(5)$.

\begin{theorem}\label{thm:triangular}
\begin{enumerate}
\item
If $l\ge4$, then $1\notin M(T(l))$.
\item
If $l$ is odd, then $2\notin M(T(l))$.
\end{enumerate}
\label{t:triang}
\end{theorem}

\begin{pf}
(a) The graph $T(l)$ has clique number $l-1$ and independence number $\lfloor l/2\rfloor$,
so cannot be embedded into its complement.

\medskip

(b) The proof is by contradiction and generalizes Schwenk's argument~\cite{schwenk} showing that three Petersen graphs cannot partition $K_{10}$. Assume that $2\in M(T(l))$ and consider
a decomposition of $2K_{{l\choose 2}}$ into
$2\left({l\choose 2}-1\right)/2(l-2)=(l+1)/2$ copies of $T(l)$.
Let $A_1,\dots,A_{(l+1)/2}$ denote the adjacency matrices of these
copies of $T(l)$. For $1\leq i\leq (l+1)/2$, denote by $\mathcal{E}_i$
the eigenspace of $A_i$ corresponding to $-2$. It is known that each
$\mathcal{E}_i$ is contained in the orthogonal complement of the all-one
vector in $\mathbb{R}^{{l\choose 2}}$ and that
$\dim(\mathcal{E}_i)={l\choose 2}-l$. Since the intersection of $m$ subspaces
each of codimension $n$ in a vector space has codimension at most $mn$, we have
\begin{eqnarray*}
\dim\left(\bigcap_{i=1}^{(l-1)/2}\mathcal{E}_i\right) 
&\geq& {l\choose 2}-1-(l-1)\cdot(l-1)/2 \\
&=& (l-3)/2 \\
&>&0.
\end{eqnarray*}
Let $x$ be a non-zero vector in
$\displaystyle{\bigcap_{i=1}^{(l-1)/2}\mathcal{E}_i}$.
Since $A_1+\cdots +A_{(l-1)/2}+A_{(l+1)/2}=2(J-I)$, we deduce that
\begin{eqnarray*}
A_{(l+1)/2}x &=& (2J-2I-A_1-\cdots-A_{(l-1)/2})x \\
&=& -2x+((l-1)/2)2x \\
&=& (l-3)x.\qquad\qquad\qquad\Box
\end{eqnarray*}
\end{pf}

\begin{prop}
$M(T(6))=4\mathbb{N}$, so that $\parti(T(6))=1$.
\end{prop}

\begin{pf}
The automorphism group $\Aut(T(6))$ is $S_6$ and has a subgroup $A_6$ which has index $7$ in the
$2$-transitive group $A_7$ of degree $15$. So $7$ copies of $T(6)$ cover the
edges of $K_{15}$ four times.\qed
\end{pf}

Thus $7$ is the smallest value of $m$ for which we don't know $M(T(m))$.
Here are a few comments on this case. Let $G=T(7)$.
\begin{enumerate}
\item By Theorem~\ref{t:triang}, we see that $1,2\notin M(G)$.
\item There does exist a $2$-$(21,6,4)$ design, namely the point residual of
the Witt design on $22$ points. However, this design is not possible in our
situation; for it has the property that any two blocks meet in $0$ or $2$
points, whereas $T(7)$ has $6$-cliques meeting in one point.
\item Applying Proposition~\ref{p:2tr}, with $H=PSL(3,4)$, we obtain
$1440\in M(G)$. This leaves a very big gap.
\end{enumerate}

\section*{Acknowledgments}
The research of the first author is partially supported by National Security Agency grant H98230-13-1-0267. The authors are grateful to Steve Butler, Darryn Bryant, Alex Fink, and Mark Walters for helpful comments which have improved this paper and to Jason Vermette for Figure 1.

\end{document}